# Decentralized formation control part I: Geometric aspects

M.-A. Belabbas *

January 14, 2011

## Abstract

In this paper, we develop new methods for the analysis of decentralized control systems and we apply them to formation control problems. The basic set-up consists of a system with multiple agents represented by the nodes of a directed graph whose edges represent an available communication channel for the agents. We address the question of whether the information flow defined by the graph is sufficient for the agents to accomplish a given task.

Formation control is concerned with problems in which agents are required to stabilize at a given distance of other agents. In this context, the graph of a formation encodes both the information flow, as described above, and the distance constraints, by fixing the lengths of the edges. A formation is said to be rigid if it cannot be continuously deformed with the distance constraints satisfied; a formation is minimally rigid if no distance constraint can be omitted without the formation losing its rigidity. Hence, the graph underlying minimally rigid formation provides just enough constraints to yield a rigid formation. An open question we will settle is whether the information flow afforded by a minimally rigid graph is sufficient to define global stabilizing control laws. We show that the answer is negative in the case of directed information flow.

In this first part, we establish basic properties of formation control in the plane. Formations and the associated control problems are defined modulo rigid transformations. This fact has strong implications on the geometry of the space of formations and on the feedback control laws, since they need to respect this invariance. We study both aspects here. In detail, we show that the space of formations of $n$ agents is $\mathbb{C}P(n-2) \times (0,\infty)$ where $\mathbb{C}P(n)$ is the complex-projective space of complex dimension $n$. We subsequently illustrate how the non-trivial topology of this space relates to the parametrization of the formation by inter-agent distances. We then establish conditions feedback control laws need to satisfy in order to yield a closed-loop system that respects both the invariance under the action of the Euclidean group and the constraints on the information flow.

---

*M.-A. Belabbas is with the School of Engineering and Applied Sciences, Harvard University, Cambridge, MA 02138 `belabbas@seas.harvard.edu`

# Contents

# 1 Introduction

Decentralized and formation control problems have been occupying a central part of the efforts in control theory for the past decade. The reason for this growing interest stems in part from the large number of potential applications—from the study of schooling and flocking, to sensor networks, to formation flight [EAM+04, SPL07]—and in part from the scientific challenges they present.

Given an ensemble of agents, a graph whose vertices are identified with the agents is used to describe the *information flow* in the system: a directed edge from agent $i$ to $j$ means that agent $i$ observes agent $j$. The main objective behind our work in formation control is to understand how the constraints on the information flow in a decentralized system affect its dynamics. Informally speaking, we would like to answer questions such as "how much does agent $i$ need to know in order for a given objective to be reached?" Here, objectives can be as varied as trying to remain in a certain part of the phase space or making sure that a given configuration is stable. We address such questions in the case of directed formation control problems. The systems involved, in addition to their intrinsic interest and in spite of their rather simple description, can have a fairly complex behavior even in low dimensions and thus provide a rich test-bed for frameworks addressing the above-mentioned issues.

In part I of this paper, we study topological and geometric properties of formations in the plane and establish a few general results about feedback control laws that respect both the invariance and the information flow of the system. In part II, we introduce an algebraic framework that captures the range of behaviors that a decentralized system can achieve and use it to study global stabilization of directed formations.

In this paper, we adopt the following point of view: cooperative control systems are most often made of a number of similar or nearly similar subsystems, and each subsystem comes with its own coordinate system. Think of a flock where each agent sees its neighbors with respect to its own position. While the parametrization of the individual systems can be made straightforward, it is often not the case for the ensemble of systems. Beyond the well-known and studied fact that interactions between systems can bring additional constraints, the parametrization of the ensemble requires more care in two aspects: the first is that the straightforward direct sum of the coordinates of each agent often results in a redundant parametrization; we will elaborate on this below. The second aspect originates from topological considerations: the dynamics of the subsystems depend on a subset of the variables necessary to describe the ensemble, as defined by the information flow in the system. We call these the localized coordinates. We will show how the interplay between these localized coordinates and the topology of the state-space of the ensemble affects the dynamics of all decentralized feedback laws for the ensemble.

Consider having, as an example, three autonomous agents in the plane called $x_1, x_2$ and $x_3$. Agent $i$ only knows the *relative position* of agent $i+1$ (taken mod 3) with respect to its own position. With these constraints, is it possible for the agents to reach an arbitrary triangular formation in the plane? It has been shown in [CMY$^+$07] that it is the case, up to the mirror symmetric of the formation, for almost all initial conditions for the agents. However, as noted in [CAM$^+$10], results beyond the ones relating to three agents have been particularly difficult to obtain.

The paper is organized as follows. In Section 2, we introduce fundamentals of graph and rigidity theory and establish most of the notation used in the paper. In particular, we introduce the mixed-adjacency and the edge-adjacency matrices that have not been used before. These operators are useful when expressing the dynamics in different coordinate systems. In Section 3, we present the necessary geometric background for the paper. It consists mainly of a review of projective spaces, the Euclidean group and the Lusternick-Schnirelmann category. The following two sections contain the main results of this paper. We establish in Section 4 geometric properties of the space of formations and relate the edge-lengths parametrization to the Lusternick-Schnirelmann category. In the last section, we introduce a model for formation control that respects the invariance properties introduced in Section 4 and derive its basic characteristics.

## 2  Frameworks and rigidity

Formation control problems are intimately related to rigidity and graph theory. We present in this section the relevant aspects and refer the reader to [GSS93] for additional details.

### 2.1  Frameworks

Let $G = (V, E)$ be a *graph* with $n$ vertices — that is $V = \{x_1, x_2, \ldots, x_n\}$ is an ordered set of vertices and $E \subset V \times V$ is a set of edges. The graph is said to be *directed* if $(i, j) \in E$ does not imply that $(j, i) \in E$. We let $|E| = m$ be the cardinality of $E$. We call the *outvalence* of a vertex the number of edges originating from this vertex and the *invalence* the number of incoming edges.

Directed graphs are used to encode the information flow in decentralized control problems. We follow the convention that an arrow leaving vertex $i$ for vertex $j$ means that agent $i$ observes agent $j$. For example, for the graph of Figure 1a, we have that $x_1$ sees both $x_2$ and $x_4$, $x_2$ and $x_4$ see $x_3$ and so forth.

The topology of a graph can conveniently be captured by some operators, instead of the sets $V$ and $E$. This is done using *adjacency matrices*:

**Definition 1.** *Given a directed graph $G = (V, E)$, its adjacency matrix $A_d \in \mathbb{R}^{n \times n}$ is defined by*

$$A_{d,ij} = \begin{cases} -1 & \text{if } (i,j) \in E \\ 0 & \text{otherwise.} \end{cases}$$

Assume that the edges are ordered. The edge-adjacency matrix is a $|E| \times |E|$ matrix whose entry $i, j$ is $-1$ if edge $i$ and edge $j$ originate from the same vertex $i$, 1 if edge $i$ ends at vertex where edge $j$ starts and 0 otherwise. Notice that $A_{e,ij}$ is zero if edge $i$ starts where edge $j$ ends and that the diagonal entries are $-1$:

**Definition 2** (Edge-adjacency matrix). *Given a directed graph $G = (V, E)$, its edge-adjacency matrix $A_e \in \mathbb{R}^{m \times n}$ is defined by*

$$A_{e,ij} = \begin{cases} -1 & \text{if } e_i = (s,t), e_j = (s,t'), s, t, t' \in V \\ 1 & \text{if } e_i = (s,t), e_j = (t,s'), s, s', t \in V \\ 0 & \text{otherwise.} \end{cases}$$

The mixed-adjacency matrix is a $|E| \times |V|$ matrix whose entry $i, j$ is $-1$ if edge $i$ originates from vertex $j$, 1 if edge $i$ ends at vertex $j$ and 0 otherwise:

**Definition 3** (Mixed adjacency matrix). *Given a directed graph $G = (V, E)$, its mixed adjacency matrix $A_m \in \mathbb{R}^{n \times m}$ is defined by*

$$A_{m,ij} = \begin{cases} -1 & \text{if } e_i = (j, s), s \in V \\ +1 & \text{if } e_i = (k, j), k \in V \\ 0 & \text{otherwise.} \end{cases}$$

The mixed- and edge-adjacency matrices are used to relate the dynamics of formation control systems in terms of the position of the agents to its expression in terms of edge lengths. We will often have to consider the matrices $A_m \otimes I$ and $A_e \otimes I$ where $\otimes$ is the Kronecker product and $I$ the two-by-two identity matrix. In order to keep the notation simple, we write $A_m^{(2)}$ and $A_e^{(2)}$ for these Kronecker products.

**Example 1.** *The mixed-adjacency matrix of the 2-cycles of Figure 1b is*

$$B = \begin{bmatrix} -1 & 1 & 0 & 0 \\ 0 & -1 & 1 & 0 \\ 1 & 0 & -1 & 0 \\ 0 & 0 & 1 & -1 \\ -1 & 0 & 0 & 1 \end{bmatrix}. \tag{1}$$

*and its edge- adjacency matrix is*

$$A_e = \begin{bmatrix} -1 & 1 & 0 & 0 & -1 \\ 0 & -1 & 1 & 0 & 0 \\ 1 & 0 & -1 & 0 & 1 \\ 0 & 0 & 1 & -1 & 0 \\ -1 & 0 & 0 & 1 & -1 \end{bmatrix} \tag{2}$$

*where edge 1 corresponds to $z_1$, etc.* □

We call a *framework* an embedding of a graph in $\mathbb{R}^2$ endowed with the usual Euclidean distance, i.e. given $G = (V, E)$ a directed graph, a framework $p$ is a mapping

$$p : V \to \mathbb{R}^2.$$

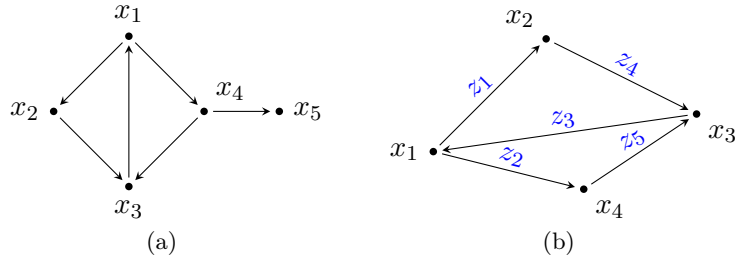

Figure 1: *(a)* A representation of the directed graph $V = \{x_1, x_2, x_3, x_4, x_5\}$, $E = \{(x_1, x_4), (x_1, x_2), (x_2, x_3), (x_3, x_1), (x_4, x_3), (x_4, x_5)\}$. *(b)* The 2-cycles formation.

By abuse of notation, we will often write $x_i$ for $p(x_i)$. We define the *distance function* $\delta$ of a framework with $n$ vertices as

$$\delta(p) : \mathbb{R}^{2n} \to \mathbb{R}_+^{n(n-1)/2} : (x_1, \ldots, x_n) \to \frac{1}{2} \begin{bmatrix} \|x_1 - x_2\|^2 \\ \|x_1 - x_3\|^2 \\ \vdots \\ \|x_1 - x_n\|^2 \\ \|x_2 - x_3\|^2 \\ \vdots \\ \|x_{n-1} - x_n\|^2 \end{bmatrix},$$

where $\mathbb{R}^+ = [0, \infty)$. We denote by $\delta(p)|_E$ the restriction of the range of $\delta$ to edges in $E$.

For a graph $G$ with $m$ edges, we define

$$\mathcal{L} = \left\{ d = (d_1, \ldots, d_m) \in \mathbb{R}_+^m \text{ for which } \exists p \text{ with } \delta(p(V))|_E = d \right\},$$

and we denote by $\mathcal{L}_0$ the interior of $\mathcal{L}$. In other words, $\mathcal{L}$ is the set of feasible assignments of edge lengths of a given graph. For example, for a triangle with edge lengths $d_1, d_2$ and $d_3$, $\mathcal{L}$ is given by the inequalities $d_1 + d_2 \geq d_3, d_2 + d_3 \geq d_1, d_1 + d_3 \geq d_2$ and $d_i \geq 0$.

When dealing with frameworks of a given graph, it is often the case that some properties are true except for a small set of frameworks, such as frameworks with all the vertices aligned, or frameworks with two vertices superposed. In order to easily deal with these particular cases, we introduce the following definition [Con05]:

**Definition 4** (Generic frameworks). *A framework is* generic *if the coordinates of all its the vertices are algebraically independent over the rationals.*

The above definition means that the positions of the vertices cannot be described as the zeros of a polynomial with rational coefficient. It is easy to see that generic frameworks are dense in the space of frameworks, and that non-generic frameworks are of measure zero in this space. We end this section with the following definition:

**Definition 5** (k-vertex connected graphs). *A graph is* connected *if for every pair of vertices $x_i, x_j \in V$, there is a path in $G$ that starts at $x_i$ and ends at $x_j$. A graph is said to be* k-vertex connected *if every graph $G'$ obtained by removing $k-1$ vertices from $V$, and the edges incident to these vertices, is connected.*

## 2.2 Rigidity

Given a graph $G$ with $m$ edges, we are given $m$ positive numbers $d_i$ and consider the frameworks of $G$ whose distance function satisfies

$$\delta(p)|_E = d_i.$$

Rigidity theory is concerned with how many edges are necessary so that a framework cannot be continuously deformed, with the exception of translations and rotations.

The *rigidity matrix* of the framework is the Jacobian $\frac{\partial \delta}{\partial x}$ restricted to the edges in $E$. We denote it by $\frac{\partial \delta}{\partial x}|_E$.

There are several notions of rigidity that are relevant for our work:

1. *Static rigidity*: A framework is said to by statically rigid, or simply rigid, if there are no nearby frameworks $p'$, modulo rotation and translation, with $\delta(p')|_E = d$.

2. *Infinitesimal rigidity*: A framework is said to be infinitesimally rigid if there are no vanishingly small motions of the vertices that keep the edge-length constraints on the framework satisfied. This translates into [GSS93]:

$$\text{rank}(\frac{\partial \delta}{\partial x}|_E) = 2n - 3.$$

3. *Minimal rigidity*: A framework is said to be minimally rigid if none of the $m$ frameworks with $m-1$ edges obtained by removing one edge is rigid.

4. *Global rigidity*: A framework is said to be globally rigid if the only other framework satisfying the same edge lengths is its mirror symmetric.

While infinitesimally rigid frameworks are rigid, the converse is *not* true, see [GSS93] for an example. The framework in Figure 1a is not statically rigid since continuous motions of $x_5$ are allowed. The framework of Figure 1b is statically rigid but not globally rigid: indeed one could take the framework where $x_2$ is sent to its mirror symmetric along $z_3$ and satisfy the same distance constraints. We will look into global rigidity in more detail in later sections.

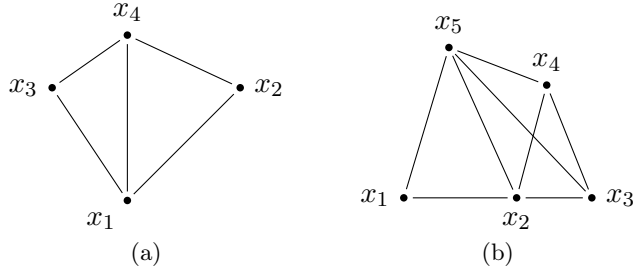

Figure 2: (a) is minimally and infinitesimally rigid, while (b) is infinitesimally rigid but not minimally rigid.

In dimension two, rigid formations are completely characterized by Laman's theorem. Given a graph with $n$ vertices and $m$ edges, a subgraph is obtained by keeping a subset of the vertices and the edges that link these vertices. Precisely, a subgraph $G'$ is given by a pair $(V' \subset V, E' \subset E)$ with $(i,j) \in E'$ if and only if $i, j \in V'$ and $(i,j) \in E$.

**Theorem 1** (Laman). *A planar framework with $n$ vertices is generically rigid if and only if*

- *it has $2n - 3$ edges*

- *for every subgraph with $n'$ vertices and $m'$ edges we have $m' \leq 2n' - 3$.*

The first part of the statement is a simple dimensionality argument: specifying $n$ agents in the plane can be done by giving their $2n$ coordinates from which we subtract 2 degrees of freedom due to invariance under translation and 1 degree of freedom due to invariance under rotation. Hence, the non-trivial part is the second part, which we understand as a density argument: the $2n - 3$ edges cannot be concentrated around too few vertices. The proof is rather involved, and we refer to [Lam70] for more information.

## 3 Projective space, Euclidean group and LS-category

We now introduce the main geometric notions that we will need in this paper. Given a framework $p$, we say that vertices are *totally coincidental* if they are all mapped to the same location by $p$:

$$p(x_1) - p(x_2) = \ldots = p(x_{n-1}) - p(x_n) = 0.$$

We consider frameworks such that the vertices are not totally coincidental for the rest of this paper, but we allow two or more vertices to be mapped to the same point. We study the set $E^n$ of configurations of $n$ non-totally coincidental agents in the plane.

## 3.1 Euclidean group

The invariance of the frameworks under rotations and translations can be formalized as an invariance under a group action. We elaborate on this here. Recall the group $SE(2)$ of affine rigid transformation of the plane, i.e. transformations consisting of a rotation and translation. Notice that we do not consider invariance under mirror symmetry at this point. This group is a three-dimensional connected Lie group. By introducing the so-called affine coordinates

$$[x_{11}, x_{12}]^T \to [x_{11}, x_{12}, 1]^T,$$

we can write a typical element of this group as

$$A(\theta, a, b) = \begin{bmatrix} \cos(\theta) & \sin(\theta) & a \\ -\sin(\theta) & \cos(\theta) & b \\ 0 & 0 & 1 \end{bmatrix}, \qquad (3)$$

with $\theta \in [0, 2\pi]$, and $a, b \in \mathbb{R}$. First, observe that

$$A(\theta, a, b) A(\theta', a', b') = A(\theta + \theta', a + a', b + b'),$$

from which we conclude that the matrices of the type of Equation (3) indeed form a group (the inverse of $A(\theta, a, b)$ is $A(-\theta, -a, -b)$).

Consider the product

$$A(\theta, a, b) \begin{bmatrix} x_{11} \\ x_{12} \\ 1 \end{bmatrix} = \begin{bmatrix} a + \cos(\theta) x_{11} + \sin(\theta) x_{12} \\ b - \sin(\theta) x_{11} + \cos(\theta) x_{12} \\ 1 \end{bmatrix};$$

the right-hand side of the above equation corresponds to the translation by $(a, b)$ of the point $x_1$ first rotated by an angle $\theta$ about the origin. In particular, the third coordinate is always one. We will write

$$A \cdot x = \begin{bmatrix} a + \cos(\theta) x_{11} + \sin(\theta) x_{12} \\ b - \sin(\theta) x_{11} + \cos(\theta) x_{12} \end{bmatrix}$$

to denote the group action of $SE(2)$ on $\mathbb{R}^2$ just described.

## 3.2 Projective spaces

In this subsection we give a brief primer on projective spaces. We recommend to the reader who is not already familiar with these spaces to read Section 4 in parallel, as the construction done in that section illustrates many of the definitions presented here. Let $V$ be a $n-$dimensional vector space; its projective space $VP(n)$ is the space of lines in $V$ passing through the origin [Mun00]. Projective spaces are examples of quotient spaces, i.e.

spaces obtained as the quotient of a manifold $M$ by a group action [Hel78]. To elaborate on this, we define $S^{n-1}$ to be the unit sphere in $\mathbb{R}^n$ centered at the origin:

$$S^{n-1} = \{u \in \mathbb{R}^n | \sum_i u_i^2 = 1\}.$$

Let $V = \mathbb{R}^n$; every line passing through the origin in $\mathbb{R}^n$ intersects $S^{n-1}$ in two antipodal points. Using this observation, we can define an action of the group $\mathbb{Z}_2$ on $S^{n-1}$ which sends a point $(x_1, \ldots, x_n)$ to its antipodal $(-x_1, \ldots, -x_n)$. There are no points in $S^{n-1}$ such that $x = -x$. This action is thus free, or without fixed point. The real projective space $\mathbb{R}P(n)$ is then defined as the set of equivalence classes of points in $S^{n-1}$ under the above-defined $\mathbb{Z}_2$ action, or similarly— since every equivalence class consists of a point in the sphere and the antipodal point— as the sphere with antipodal points identified:

$$\mathbb{R}P(n) = S^{n-1}/\mathbb{Z}_2 = \left\{x \in S^{n-1} \text{ with } x \sim -x\right\},$$

where $a \sim b$ means that $a$ and $b$ are identified.

Projective spaces appear most commonly in control theory and signal processing as instances of Grassmanian manifolds, and there is large body of literature dealing with how to perform numerical computations and integration on such spaces. They are rather difficult to visualize, but there is a useful recursive construction for them which we describe here.

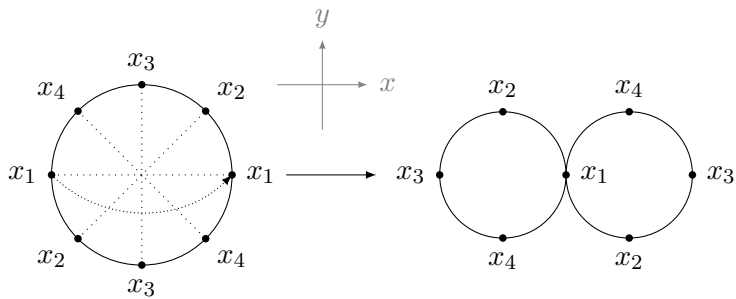

Figure 3: An action of $\mathbb{Z}_2$ send a point $x$ on the circle to its antipodal $-x$. The quotient $S^1/\mathbb{Z}_2 = \mathbb{R}P(1)$ is the set of equivalence classes of points under this $\mathbb{Z}_2$ action. If we identify the two elements in the equivalence class of $x_1$, we obtain two identical circles tangent at $x_1$. These circles contain the same points and each is thus a copy of $\mathbb{R}P(1)$

Let us first consider the space $\mathbb{R}P(1)$ of lines in $\mathbb{R}^2$. By the previous paragraphs, it is obtained by identifying antipodal points of the circle $x^2 + y^2 = 1$ in $\mathbb{R}^2$. We can represent every equivalence class by a point in the upper half plane (i.e. by its representative with $y > 0$.), with the exception of the equivalence class of points with $y = 0$, where the above rule does not yield a unique choice; see Figure 3. Identifying the points $x_1$ and $-x_1$ turns the circle into a figure 8, with its two circles tangent at $x_1$. These are thus two copies of $\mathbb{RP}(1)$, since they both contain the same points. We conclude that

$$\mathbb{RP}(1) \simeq S^1.$$

While the identification of projective spaces with spheres does not generalize to higher dimensions, the construction does.

### 3.2.1 $\mathbb{RP}(n)$ for $n \geq 2$

Using the same approach for $\mathbb{RP}(2)$, we start with the sphere $x_1^2+x_2^2+x_3^2 = 1$ in $\mathbb{R}^3$ and $\mathbb{RP}(2)$ is the space of equivalence classes under the $\mathbb{Z}_2$ action that sends a point to its antipodal. We can similarly choose to represent equivalence classes by their representative with $x_3 > 0$—this space is a disk of dimension two, or $D^2$— which yields a unique representative except for points on the circle $x_3 = 0$, since both points in the equivalence class have the third coordinate zero. But by the previous section, the identification of antipodal points on the circle yields $\mathbb{RP}(1)$, and then we can represent $\mathbb{RP}(2)$ as a disk $D^2$ with its boundary $S^1$ having its antipodal point identified to yield $\mathbb{RP}(1)$. This construction, which is illustrated in Figure 4, generalizes in a straightforward manner to higher dimensions: $\mathbb{RP}(n)$ can be seen as a ball $D^n$ with an $\mathbb{RP}(n-1)$ boundary.

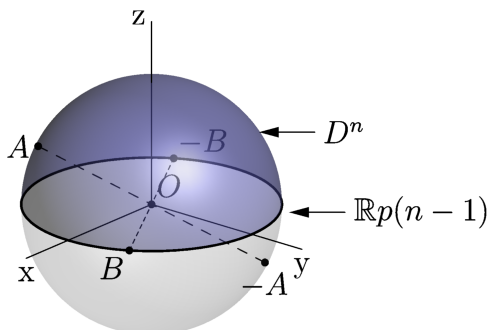

Figure 4: The action of $\mathbb{Z}_2$ on the sphere sends a point $x$ to its antipodal $-x$. The quotient $S^n/\mathbb{Z}_2$ is the space of equivalence classes under this action. We keep the representatives of these equivalence classes with $x_3 > 0$, e.g. $A$ over $-A$. They are in the upper half-sphere. For points on the equator, e.g. $B$ and $-B$, both representatives have $x_3 = 0$. Identifying these point corresponds to defining a copy of $\mathbb{RP}(n-1)$ at the equator.

### 3.2.2 Complex projective spaces

Complex projective spaces are defined similarly to real projective spaces, but the underlying vector spaces are complex. Hence, we are looking at the space of complex lines (i.e. copies of $\mathbb{C}$) through the origin in the complex vector space $\mathbb{C}^n$. We can express $\mathbb{C}P(n)$ as a quotient space as follows: given[1] $z \in \mathbb{C}_0^{n+1}$ with $z = (z_0, \ldots, z_n), z_i \in \mathbb{C}$, the complex line passing through 0 and $z$ is the set of points $z' = (\lambda z_0, \ldots, \lambda z_n)$ for $\lambda \in \mathbb{C}_0$. Thus $\mathbb{C}P(n)$ is given by the quotient

$$\mathbb{C}P(n) \simeq \mathbb{C}_0^{n+1}/\mathbb{C}_0$$

where the action of $\mathbb{C}_0$ is the multiplication by $\lambda$ given above. An element in $\mathbb{C}P(n)$ is often given by a point in $\mathbb{C}_0^{n+1}$, with the understanding that this point is defined up to multiplication by $\lambda \in \mathbb{C}_0$; these are the so-called *homogeneous coordinates* and they are denoted by a bracket notation $[z_0 : \ldots : z_n]$ to indicate that $(z_0, z_1, \ldots, z_n)$ is a representative of an equivalence class:

$$[z_0 : \ldots : z_n] \simeq [\lambda z_0 : \ldots : \lambda z_n], \ \forall \lambda \in \mathbb{C}_0.$$

Observe that having all the $z_i = 0$ is not allowed since in that case $\lambda z$ does not define a line. Using homogeneous coordinates, we can exhibit copies of $\mathbb{C}P(n-1)$ in $\mathbb{C}P(n)$, similarly to what was done in the real case. First, observe that if $z_0 \neq 0$, we can divide the homogeneous coordinates by $z_0$ to obtain $[1 : z_1/z_0 : \ldots : z_n/z_0]$, and if $z_0 \neq 0$, all the $z_i, i \neq 0$ can be zero simultaneously since $z_0 = 1$ ensures that not all coordinates are zero. Thus if $z_0 \neq 0$, the coordinates $z_1, \ldots z_n$ are in $\mathbb{C}^n$ without restriction. This space is topologically equivalent to a disk $D^{2n}$ of real dimension $2n$. When $z_0 = 0$, the homogeneous coordinates are of the form $[0 : z_1 : \ldots : z_n]$: they are the coordinates of a copy of $\mathbb{C}P(n-1)$. Hence, we have a decomposition, similar to the one of $\mathbb{R}P(n)$, of $\mathbb{C}P(n)$ into a disk $D^{2n}$ and with $\mathbb{C}P(n-1)$ at its boundary.

Using the above construction recursively, we obtain

$$\mathbb{C}P(n) = D^{2n} \cup D^{2n-2} \cup \ldots \cup D^0$$

where $D^0$ is a point. Each time a disk is added in the construction, some identifications of points, as described in the above paragraph, have to be made.

## 3.3 The Lusternick-Schnirelmann category

The Lusternich-Schnirelmann category of a manifold, or LS-category, is an important topological invariant that is related to the minimal number of charts contained in any atlas for the manifold. We describe it in this section.

A set $S$ in $\mathbb{R}^n$ is homotopic [Mun00] to a point if there exists a continuous map $F(t, x) : [0, 1] \times \mathbb{R}^n \to \mathbb{R}^n$ such that $F(0, x)$ is the identity on $S$ and $F(1, x)$ maps $S$ to a point. Let $M$ be a smooth manifold of dimension $n$ and $\mathcal{V} = \{V_i, i = 1 \ldots k\}$ be a collection of closed connected sets which are *homotopic to a point*, or null-homotopic. Convex sets, or sets that are homeomorphic to a convex set, are homotopic to a point; sets that are not null-homotopic include sets which have "holes" in them. In dimension 1, the null-homotopic sets are the connected intervals. In $\mathbb{R}^2$, null-homotopic sets are the *simply connected sets*: sets such that every closed path contained in them can be continuously deformed to a point (see Figure 5).

---

[1] $\mathbb{C}_0^n = \mathbb{C}^n - \{0\}$

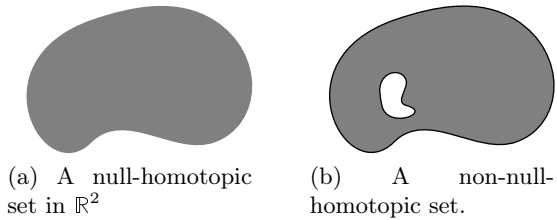

(a) A null-homotopic set in $\mathbb{R}^2$

(b) A non-null-homotopic set.

Figure 5

Let $I$ be an index set and

$$\mathcal{V}_I = \{V_i, i \in I | V_i \text{ is a closed set}\}$$

a collection of closed sets indexed by $I$. We say that the collection $\mathcal{V}_I$ is a cover for $M$ if $M \subset \cup_{i \in I} V_i$. The Lusternick-Schnirelmann category [CLOT03] of $M$, or LS-category, is the topological invariant defined as follows:

**Definition 6** (Lusternick-Schnirelmann Category). *Given a closed manifold $M$, the Lusternick-Schnirelmann category of $M$, written $\text{cat}(M)$, is the the least cardinality $|I|$, over all possible $\mathcal{V}_I$ that are covers for $M$.*

For example, it is easy to see that one needs at least two connected intervals to cover the circle $S^1$. Hence $\text{cat}(S^1) = 2$. Actually, one can prove that $\text{cat}(S^n) = 2$: two disks of dimension $n$ are needed to cover $S^n$. Even though the original definition of the LS-category was in terms of closed covers as introduced above, it is nowadays more common to encounter a definition in terms of open covers, i.e. covers where the sets $V_i$ are open. We refer to [CLOT03] for more information and relations between the two quantities.

For our purposes, we will need the Lusternick-Schnirelmann category of complex projective spaces, which is known to be [CLOT03]

$$\text{cat}(\mathbb{C}P(n)) = n + 1.$$

## 4 The Geometry of the space of formations

We define $E^n$ to be the space of equivalence classes, under rotations and translations, of formations of $n$ agents in the plane. We use the notions introduced in the previous section to characterize these spaces. For a different point of view with an eye on statistical applications, we refer to the excellent survey [Ken89]. We start by the space $E^3$ of three agents in the plane. We work out this case in details as it sheds light on the more abstract constructions needed for $n$ agents. We recall here that the agents in the formation are *labelled*.

## 4.1 The space $E^3$

Consider three agents $x_1$, $x_2$ and $x_3$ in the plane $\mathbb{R}^2$. We can describe their position with a vector $x \in \mathbb{R}^6$ with the first two coordinates of $x$ giving the position of $x_1$, the next two coordinates the position of $x_2$ and so forth. We will write for the rest of this section $x = [x_1, x_2, x_3]^T$ for $x \in \mathbb{R}^6$. Since formations are defined up to a translation and rotation, it is clear that the above coordinates on the set of formations with 3 agents is redundant. This redundancy takes the form of a group action of $SE(2)$ on $\mathbb{R}^6$ given by

$$A \cdot x = \begin{bmatrix} Ax_1 \\ Ax_2 \\ Ax_3 \end{bmatrix} \quad (4)$$

where the action of $A \in SE(2)$ on $x_i \in \mathbb{R}^2$ is the one given in the previous section.

The space of totally coincidental formations is the subspace of $\mathbb{R}^6$ defined as

$$\mathcal{N} = \left\{ x \in \mathbb{R}^6 | x_1 = x_2 = x_3 \right\}.$$

*Fact 1*: $\mathcal{N}$ is invariant under the action of $SE(2)$.
This fact is obvious since $x_i = x_j$ implies $Ax_i = Ax_j$.

*Fact 2*: $E^3$ is connected.
Since $\mathcal{N}$ is a linear subspace of codimension four, it does not separate $\mathbb{R}^6$ in two regions and $E^3$ has a single connected component.

We now show that the space $E^3$ is actually isomorphic to $\mathbb{R}^3 - \{0\}$, but the mapping that sends a point in $\mathbb{R}^3$ to a triangle is not obvious. We will later give a more abstract proof of this fact.

**Proposition 1.** *The space $E^3$ of triangular formations in the plane is $\mathbb{R}^3_0$.*

*Proof.* We work out explicitly the quotient $(\mathbb{R}^6 - \mathcal{N})/SE(2)$. Given an equivalence class of triangles, we can always choose a representative with $x_1 = (0,0)$, which takes care of the translational redundancy. If we assume that $x_1 \neq x_2$, we can take a representative such that the $x_1 - x_2$ edge is aligned with the $x-$ axis in $\mathbb{R}^2$. Hence, we can describe a non-degenerate triangle by a point $u \in \mathbb{R}^3_0$ with $x_1 = (0,0), x_{22} = 0$ and $x_{21} = u_1$, $x_{31} = u_2$ and $x_{32} = u_3$. This representation is still redundant: a point $u$ yields a triangle which is equivalent by a rotation of 180 degrees to the formation given by $-u$. The quotient of $\mathbb{R}^3_0$ by $u \simeq -u$, i.e. the space of formations with $x_1$ not equal to $x_2$, is $\mathbb{RP}(2) \times \mathbb{R}^+$ as described in the previous section.

We now show how the degenerate formations are attached to the space of non-degenerate formation just described. If $x_1 = x_2$, then $x_2 \neq x_3$—otherwise the formation is totally coincidental. Assume we have $u_1 = 0$ and $(u_2, u_3) \neq (0,0)$. Geometrically, this corresponds to choosing the point $x_3$ in the plane minus the origin. This point can be described in polar coordinates as $x_3 = re^{i\theta}$.

To simplify matters, we assume that the formation is normalized in the sense that the sum of the edge lengths is fixed. It is easy to see that this fixes a value in $\mathbb{R}^+$ for coordinates in the product $\mathbb{RP}(2) \times \mathbb{R}^+$. Hence, the space of (normalized) formations with $x_2 \neq x_1$ is $\mathbb{RP}(2)$ and the formations with $x_2 = x_1$, which are described by $z_3 = re^{i\theta}$, are a circle on its boundary. The situation is similar to the one depicted in Figure 4, where the upper half-sphere corresponds to formations with $x_2 \neq x_1$ and the equator to formations with $x_2 = x_1$. Since the formations are defined up to a rotation, all formations with $x_2 = x_1$ are identified. This corresponds to shrinking the circle in Figure 4 to a point, which yields two identical spheres $S^2$ touching at a point similarly to the situation of Figure 3. Repeating this construction for every normalization coefficient in $\mathbb{R}^+$, we conclude that the space of formations is $S^2 \times \mathbb{R}^+ \simeq \mathbb{R}^3_0$. ∎

## 4.2 Parametrization of $E^3$ by edge lengths

We now look at the relation between the localized coordinates $\mathcal{L}$ and the geometry of $E^3$. Let us write $d_1 = \|x_2 - x_1\|, d_2 = \|x_3 - x_2\|$ and $d_3 = \|x_3 - x_1\|$. Using the result of Proposition 1, we know we can consider the space of formations where, say, the sum of edge lengths $d_1 + d_2 + d_3 = 1$, and the space of such formations is then the 2−sphere $S^2$.

The edge lengths have to satisfy the triangle inequalities

$$d_1 \leq d_2 + d_3, d_2 \leq d_3 + d_1, d_3 \leq d_1 + d_2. \tag{5}$$

Thus the set of admissible $d_i$ forms a convex subset of $\mathbb{R}^3$, whose edges are such that the above inequalities are equalities (see Figure 6).

**Proposition 2.** *There is an open set $U \subset \mathcal{L}$ such that there are at least two frameworks that correspond to $(d_1, d_2, d_3) \in U$.*

The proof establishes a link between the number of non-congruent formations and the topology of $S^2$. We will generalize this result below.

*Proof.* The set described by the inequalities (5) is a cone in $\mathbb{R}^3$. It is easy to see that an intersection of this cone with a plane orthogonal to the main diagonal[2] is a convex set which corresponds to having the sum of the edge lengths normalized. Precisely, if the orthogonal plane intersects the main diagonal at $(a, a, a)$, then $d_1 + d_2 + d_3 = a$; we call this intersection $\mathcal{L}_a$. Fixing the value of $a$ corresponds to fixing the value of the variable corresponding to $\mathbb{R}^+$ in the description of $E^3$. Without loss of generality, we consider the space of *normalized* formations in three agents, which is $S^2$ by Proposition 1.

---
[2]we call *main diagonal* in $\mathbb{R}^n$ the subspace spanned by the vector $(1, 1, \ldots, 1)$

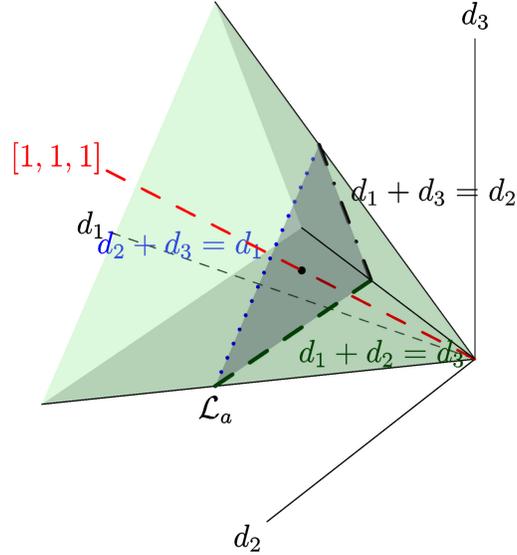

Figure 6: The set $\mathcal{L}$ for triangular formations is a cone in $\mathbb{R}_+^3$. We depict its intersection with a plane with normal $[1, 1, 1]$ at $[a, a, a]$. We call this intersection $\mathcal{L}_a$.

Consider the map
$$\nu : S^2 \to \mathcal{L}_a : p(V) \to \delta(p(V))$$
that maps frameworks to edge lengths. Because $\nu$ is continuous, it is enough to show that it is not injective to prove the existence of $U$.

By definition of $\mathcal{L}$, $\nu$ is onto. Reasoning by a contrapositive argument, assume that $\nu$ is injective; in other words, there is only one framework for every given vector $(d_1, d_2, d_3) \in \mathcal{L}_a$. This means that we can uniquely assign a point in $S^2$ to every $(d_1, d_2, d_3)$ and vice-versa. But this is in contradiction with the fact that the Lusternick-Schnirelmann category of $S^2$ is 2: closed convex sets are null-homotopic and the above claim is equivalent to covering $S^2$ with a single such set. ∎

It is easy to see that in the case of triangles, there are *exactly* two such formations; they correspond to a triangle and its mirror symmetric, as illustrated in Figure 7.

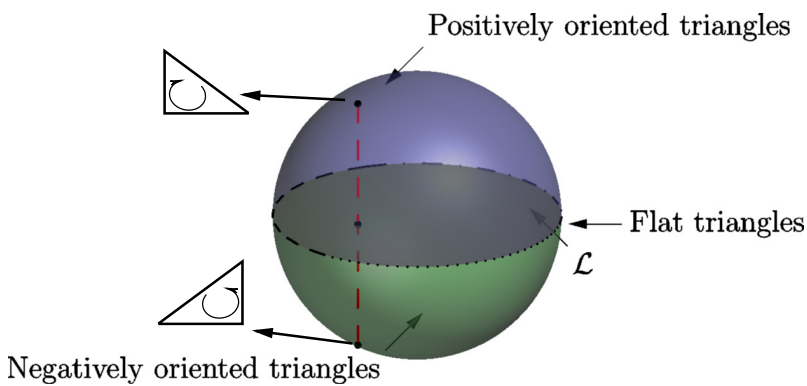

Figure 7: The space of normalized triangles in the plane is $S^2$. The set of normalized edge lengths corresponds to a triangle, depicted in Figure 6, which we have deformed above into a disk $\mathcal{L}_a$. To each point in the disk correspond two triangles with the same edge lengths but different orientations. The map $\nu$ of Proposition 4.2 is a projection of the sphere onto the disk $\mathcal{L}_a$. The knowledge of the localized coordinates (a point in the disk) is thus not sufficient to know the configuration of the ensemble (a point in $S^2$).

## 4.3 The space $E^n$

We describe the space of configurations of $n-$agents in the plane.

**Theorem 2.** *The space of equivalence classes of formations of $n$ agents in the plane $E^n$ is*

$$E^n \simeq \mathbb{C}\mathrm{P}(n-2) \times \mathbb{R}^+.$$

*Proof.* Recall that a general formation in the plane has $2n-3$ degrees of freedom, two for each of the $n$ agents and three degrees of freedom are lost to translational and rotational invariance. The projective space $\mathbb{C}\mathrm{P}(n-2)$ has a complex dimension of $n-2$ and real dimension of $2n-4$, and thus $\dim(\mathbb{C}\mathrm{P}(n-2) \times \mathbb{R}^+) = 2n-3$.

Without loss of generality, we identity the plane $\mathbb{R}^2$ with $\mathbb{C}$, hence the positions of the agents are given by $n$ complex numbers $z_1 = x_{11} + jx_{12}, \ldots, z_n = x_{n1} + jx_{n2}$. In complex coordinates, the rotation of a formation by an angle $\theta$ corresponds to multiplying the coordinates $z_i$ by $e^{j\theta}$.

Similarly to what was done in the case of $n = 3$ agents, we use the translational invariance to put the first agent at $0 + j0$. Thus a formation is described by $n-1$ complex numbers $(z_2, z_3, \ldots, z_n)$ with the identification

$$(z_2, z_3, \ldots, z_n) \simeq e^{j\theta}(z_2, z_3, \ldots, z_n), \forall \theta.$$

Recall that the set obtained via the identification $(z_2, z_3, \ldots, z_n) \simeq re^{j\theta}(z_2, z_3, \ldots, z_n), \forall \theta \in [0, 2\pi], r \in \mathbb{R}^+$ is $\mathbb{C}\mathrm{P}(n-2)$. Hence, we have

$$E^n = \mathbb{C}\mathrm{P}(n-2) \times \mathbb{R}^+$$

as announced. ∎

## 4.4 Parametrization by edge lengths

Formation control is concerned with finding feedback laws that stabilize agents at a given distance from a subset of the other agents, doing so based only on partial information about the formation. It is thus of interest to know how many formations satisfy the distance constraints specified.

The first requirement that was encountered was the one of rigidity: if the formation is not rigid, there exists a continuous set of formations that satisfy the specified edge lengths. The study of these manifolds has a long history that can be traced back to at least Lippman Lipkin and Peaucelier. They were investigating the design of non-rigid frameworks that would transform a circular motion to a straight-line motion. In this context, the framework is understood as an input/output system, and some of the vertices have fixed positions. More recently, William Thurston proved that using a complex enough framework, any curve in the plane can be approximately obtained as an output of a non-rigid framework [Thu97], or that one "could sign one's name" with a planar framework.

Recall that if a framework is rigid, it does *not* in general imply that only this framework and its mirror symmetric satisfy the given edge length constraints. We called frameworks for which this is true *globally rigid*. Before discussing this point any further, we formally define the mirror symmetry of a framework $x$ to be the operation that sends

$$x_i = \begin{bmatrix} x_{i1} \\ x_{i2} \end{bmatrix} \to \begin{bmatrix} x_{i1} \\ -x_{i2} \end{bmatrix}.$$

Notice that the mirror symmetric of a framework is not equivalent to the original formation via a rigid transformation since it reverses the orientation.

Unlike minimal rigidity, global rigidity is a function of the framework and not of the underlying graph $G$. Moreover, deciding whether a given framework is globally rigid is an NP-hard problem [EBM79], even for one dimensional frameworks. In order to avoid complications inherent to this dependence on a particular framework, one considers *generically globally rigid* formations, as they can be characterized in terms of their underlying graphs:

**Theorem 3** (Hendrickson, Connelly). *A graph $G$ with $n \geq 4$ vertices is* generically globally rigid *in $\mathbb{R}^2$ if and only if $G$ is generically redundantly rigid and vertex 3-connected.*

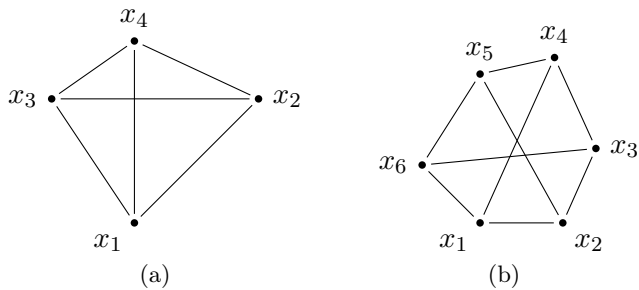

Figure 8: (a): A generically globally rigid framework. (b): A globally rigid framework that is not generically globally rigid [Con05]. It is globally rigid if and only if the vertices lie on a conic in $\mathbb{R}^2$. Observe that it is also minimally rigid.

A graph is generically redundantly rigid if, for any generic framework with underlying graph $G$, deleting one edge leaves the framework rigid. Hence, the *minimally rigid* formations, which are the main object of study in cooperative control, *are not generically globally rigid*. This yields a lower bound of 4 on the number of equilibria of formation control system with $n \geq 4$ agents solely relying on the inter-agent distance to stabilize a minimally rigid formation. Indeed, the mirror symmetry insures that there always is an even number of frameworks with a given edge lengths vector, and this number is two only for globally rigid formations.

A finer bound on the number of generic non-congruent formations in the plane for a given graph can be obtained by relating topological characteristics of $\mathcal{L}$ and $\mathbb{CP}(n)$, as we show below.

**Theorem 4.** *Let $G$ be a rigid graph on $n \geq 4$ vertices. If the set of localized coordinates $\mathcal{L}$ is null-homotopic, then there are at least $2\lceil \frac{(n-1)}{2} \rceil$ frameworks for a generic edge length vector and, in particular, $G$ is not globally rigid.*

Recall that $\lceil \alpha \rceil$ is the smallest integer with $\alpha \leq \lceil \alpha \rceil$. This result gives a lower bound on the number of non-congruent formations and an obstruction to global rigidity: *a formation cannot be globally rigid if $\mathcal{L}$ is null-homotopic.*

*Proof.* We first observe that for any feasible $d \in \mathcal{L} \subset \mathbb{R}_+^m$ and $\alpha \in \mathbb{R}_+$, we have $\alpha d \in \mathcal{L}$. We can thus normalize the sum of the $d_i$ by taking the intersection of $\mathcal{L}$ with the hyperplane orthogonal to the unit diagonal in $\mathbb{R}^m$ and passing through $\alpha > 0$, similarly to what was done in previous sections. We call this intersection $\mathcal{L}_\alpha$. We have that $\mathcal{L}_\alpha$ is non-empty if $\mathcal{L}$ is non-empty and, by definition, $\sum_{i=1}^m d_i = \alpha$ for $d \in \mathcal{L}_\alpha$.

We now prove that $\mathcal{L}_\alpha$ is closed. Consider the map

$$\Phi : \mathcal{L}_\alpha \to \mathbb{CP}(n-2)$$

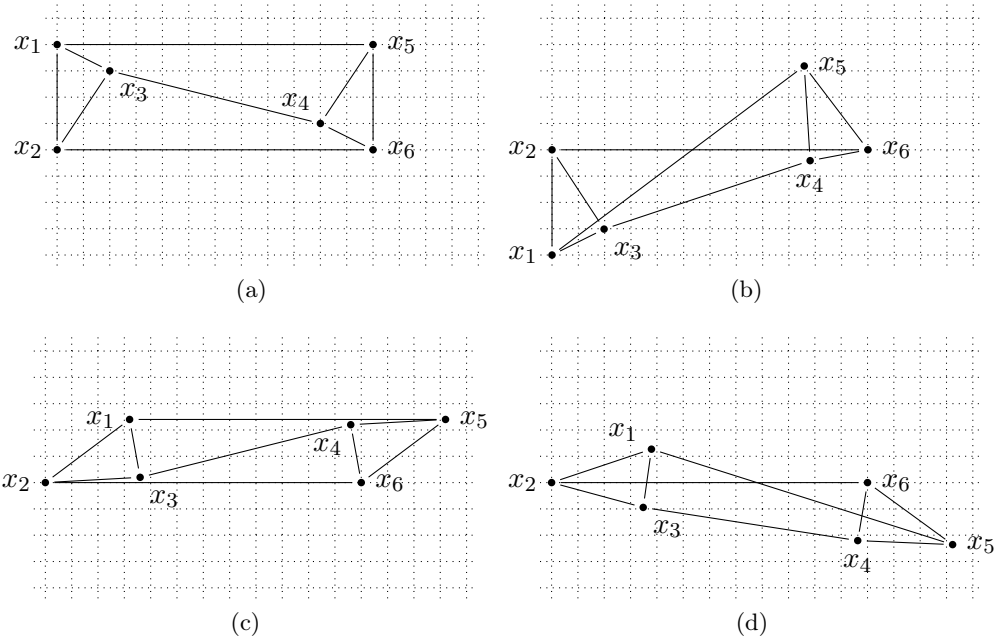

Figure 9: Four non-congruent frameworks underlying the same minimally rigid graph on 6 vertices in $\mathbb{R}^2$.

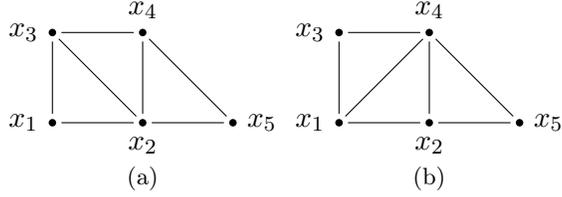

Figure 10: The two frameworks are minimally rigid in the plane, but are not isomorphic as graphs since $x_3$ has a different degree in (a) and (b).

that takes $d$ to a framework in the plane with edge lengths given by $d$. Because $G$ is rigid, there is a finite number of such frameworks. The map $\Phi$ is multi-valued and continuous on every branch of the image. This map is also onto: indeed, for every point in $\mathbb{CP}(n-2)$, we can pick a formation in the plane as in Theorem 2 with an arbitrary orientation. The pre-image of this point under $\Phi$ is obtained by reading the relevant edge lengths. Because $\Phi$ is continuous, the preimage of a closed set under $\Phi$ is closed, and thus $\mathcal{L}_\alpha$ is closed.

Additionally, if $\mathcal{L}$ is null-homotopic, then clearly $\mathcal{L}_\alpha$ is null-homotopic for all $\alpha$. Hence, if for a generic $d \in \mathcal{L}_\alpha$ there are $k$ distinct formations in the plane, then we have a closed cover of $\mathbb{CP}(n-2)$ by $k$ closed sets. Recall that the LS-category is a lower bound on the cardinality of closed covers on a space and thus $\text{cat}(\mathbb{CP}(n-2)) = n-1 \Rightarrow k \geq (n-1)$. Because the number of frameworks for a given $d$ is always even, we have that $k \geq 2\lceil \frac{n-1}{2} \rceil$. ∎

## 4.5 The space $E^4$

We explicitly describe the space $E^4$ of frameworks with four vertices in the plane. Observe that all minimally rigid frameworks with 4 agents are of the type of the (undirected) 2-cycles depicted in Figure 1b. This uniqueness property is lost for minimally rigid frameworks with $n \geq 5$ vertices, as one can easily exhibit non-equivalent minimally rigid frameworks with $n = 5$ (see Figure 10).

The constraints on the edge lengths for the 2-cycles are given by

$$\begin{aligned} d_3 &\leq d_1 + d_2 \\ d_1 &\leq d_2 + d_3 \\ d_2 &\leq d_3 + d_1 \\ d_3 &\leq d_4 + d_5 \\ d_4 &\leq d_5 + d_3 \\ d_5 &\leq d_3 + d_4 \\ d_i &\geq 0 \end{aligned}$$

From these relations, we see that $\mathcal{L}$ is a convex set. In fact, it is a cone in the positive orthant of $\mathbb{R}^5$. Summing all the inequalities above, we have $\sum_{i=1}^{5} d_i + d_3 \leq 2\sum_{i=1}^{5} d_i + 2d_3$, which simplifies to $\sum_{i=1}^{5} d_i \geq 0$. Hence, there are only five independent inequalities. Moreover, some become equalities on the hyperplanes that define the positive orthant in $\mathbb{R}^5$: for example, setting $d_1 = 0$, the first and third inequalities yield $d_2 = d_3$, etc. We conclude that the picture is similar to the one for the three agent frameworks of Figure 6, albeit in dimension five.

According to Theorem 4, there are at least $2\lceil \frac{4-1}{2} \rceil = 4$ frameworks in the plane satisfying a generic edge lengths vector $d \in \mathcal{L}$. By inspection, it is easy to see that there are exactly four such frameworks. It is informative to relate these frameworks to a discrete group action on $E^4$.

Given a framework in the plane with vertices $x_1, \ldots x_4$ and such that $d_3 > 0$, we define:

$$z_1 = x_2 - x_1, z_2 = x_3 - x_2, z_3 = x_1 - x_3, z_4 = x_3 - x_3, z_5 = x_4 - x_1, \tag{6}$$

as in Figure 1b, and let $z_3^\perp$ be the orthogonal unit vector to $z_3$. We define the *reflection symmetry* along $z_3$ as

$$\mathcal{R}^{z_3}(x) = x - 2\langle x, z_3^\perp \rangle z_3^\perp, x \in \mathbb{R}^2,$$

where $\langle \cdot, \cdot \rangle$ is the usual inner product in $\mathbb{R}^2$. We then define the two operations

$$\begin{aligned} R_1(x_1, x_2, x_3, x_4) &= (x_1, \mathcal{R}^{z_3}(x_2), x_3, x_4) \\ R_2(x_1, x_2, x_3, x_4) &= (x_1, x_2, x_3, \mathcal{R}^{z_3} x_4). \end{aligned}$$

We clearly have that

$$R_i^2 = 1 \text{ and } R_1 R_2 = R_2 R_1. \tag{7}$$

Consider the discrete group $\mathbb{Z}_2 \times \mathbb{Z}_2$ with elements $(1,0), (0,1), (1,1), (0,0)$ and the group operation is addition modulo 2. We can identify $R_1$ and $R_2$ with the elements $(1,0)$ and $(0,1)$: if we understand composition of symmetries as addition of elements in $\mathbb{Z}_2$, then Equation (7) tells us that the symmetries $R_1$ and $R_2$ obey the same multiplication table as $(1,0)$ and $(0,1)$. In this context, the identity for frameworks corresponds to the additive identity $(0,0) \in \mathbb{Z}_2 \times \mathbb{Z}_2$. Indeed $R_1 R_1 = 1 \leftrightarrow (1,0)+(1,0) = (0,0)$, $R_1 R_2 \leftrightarrow (1,0)+(0,1) = (0,1)+(1,0) \leftrightarrow R_2 R_1$, etc. We also observe that the elements $(0,0), (1,1) \in \mathbb{Z}_2 \times \mathbb{Z}_2$ form a proper subgroup which corresponds geometrically to the mirror symmetric of the formation. This is summarized in the table below, where $I$ is the identity:

|       | $R_1 = (1,0)$     | $R_2 = (0,1)$     |
|-------|-------------------|-------------------|
| $R_1$ | I=(0,0)           | $R_1 R_2 = (1,1)$ |
| $R_2$ | $R_1 R_2 = (1,1)$ | I=(0,0)           |

The section of the cone $\mathcal{L}$ obtained by fixing the sum of the edge lengths, similarly to Figure 6, is a simplex of dimension 4. This simplex is a closed null-homotopic set, and four copies of it are needed to cover $\mathbb{C}P(2)$ as was shown in Theorem 4. We represent the four 4-simplices that cover $\mathbb{C}P(2)$ in Figure 11. This figure is the equivalent of Figure 7 for the case of $n = 4$ agents.

Each of the 4 rectangles in the figure corresponds to a 4-simplex and contains frameworks related to frameworks in the other 4-simplices via the symmetries $R_1$ and $R_2$. We call $LR$ the simplex corresponding to frameworks where 2 is on the left of the 1-3 axis and 4 on the right, $LL$ the simplex where both 2 and 4 are on the left, etc. We illustrate how the simplices intersect on some codimension 1 faces, corresponding to the invariant space of $R_1$ and $R_2$—i.e. the space of frameworks which are their own symmetric under $R_1$ or $R_2$. The intersection of these invariant spaces is a codimension 2 space which corresponds to frameworks invariant under both $R_1$ and $R_2$. This invariance is tantamount to invariance under mirror symmetry. Since the sum of the edge lengths has been normalized, each $d_i$ has an upper bound. Observe that by taking the limit as $d_2$ grows to its upper bound, we find that the left boundary of the figure contains frameworks which have three vertices coincidental. But these frameworks are also invariant under $R_1$ and $R_2$ and thus there is an identification between parts of this boundary and the center corresponding to the frameworks invariant under mirror symmetry. In general, since we are representing a five dimensional space in two dimensions, we have to make some choices as to which characteristics of the space to represent. In Figure 11, we will not describe the boundaries into detail, except to mention that they contain degenerate frameworks and that intricate identifications happen there.

## 4.6 Henneberg sequences and group actions

In this section, we elaborate on a possible relation between Henneberg sequences, defined below, and discrete group actions on minimally rigid frameworks.

We have shown that in the case $n = 4$, the four frameworks that satisfy a given set of edge lengths can be characterized as the orbit of the group $\mathbb{Z}_2 \times \mathbb{Z}_2$ acting on any formation that satisfies these edge lengths. We did so by assigning to each generator of the group a reflection in the plane (namely, $R_1$ was associated with $(1, 0)$ and $R_2$ with $(0, 1)$ ) and verified that these two indeed obeyed the same multiplication table as $\mathbb{Z}_2 \times \mathbb{Z}_2$.

A major result about minimally rigid formations in the plane is the Henneberg Theorem and the associated *Henneberg* sequence [GSS93]. They state that all minimally rigid framework in the plane can be obtained inductively, starting with a single line segment (two vertices and one edge) and applying at each time step one of the two operations described below. Given a framework at step $n-1$ with a set of vertices $V_{n-1}$ and a set of edges $E_{n-1}$ of cardinalities $n-1$ and $m-1$ respectively, we define:

- *vertex add:* add a vertex to the framework and link it to two distinct existing vertices. The choice of vertices is arbitrary. Specifically, for $i \neq j$, $i, j \leq n-1$, we have:

$$\begin{aligned} V_n &= V_{n-1} \cup \{x_n\} \\ E_n &= E_{n-1} \cup \{(x_n, x_i), (x_n, x_j)\}. \end{aligned}$$

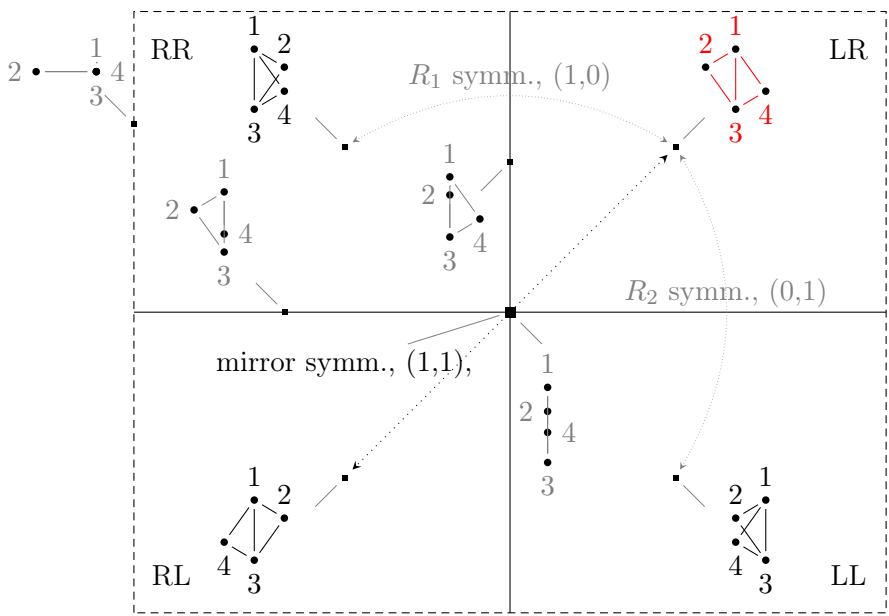

Figure 11: A depiction of the cover of $\mathbb{C}P(2)$ by four 4-simplices, depicted here as rectangles. Each simplex contains frameworks which can be identified by the position of 2 and 4 with respect to the 1-3 axis. For example, the top-right simplex contains frameworks with 2 on the left and 4 on the right. The vertical axis corresponds to a facet (i.e. subspace of codimension 1) of the simplices which contains frameworks that are invariant under $R_1$. Similarly, the horizontal axis corresponds to frameworks that are invariant under $R_2$. The intersection is a subspace of codimension 2 which corresponds to formations invariant under mirror symmetry.

- *Edge-split:* choose an edge in $E_{n-1}$ and delete it from the framework. Add a vertex and link it to the two vertices to which the deleted edge was incident and to a third, distinct vertex. Specifically, for $(x_i, x_j) \in E_{n-1}$ and $k \neq i, j$, we have:

$$V_n = V_{n-1} \cup \{x_n\}$$
$$E_n = E_{n-1} \cup \{(x_n, x_i), (x_n, x_j), (x_n, x_k)\} - \{(x_i, x_j)\}.$$

We illustrate the operations in Figure 12. The sequence of operations obtained is called *Henneberg sequence*. It has the property that all frameworks in the sequence are minimally rigid. There is no unique Henneberg sequence leading to a minimally rigid framework in general.

Consider the Henneberg sequence leading to the 2-cycles that is depicted in Figure 13 The first operation is a vertex-add applied to a single segment. This operation can be realized in two different ways: either the vertex is added to the left of the segment or to its right. This choice results in the group of symmetry $\mathbb{Z}_2$ acting on the resulting triangular framework. There is an identical choice for the second vertex-add operation which yields $\mathbb{Z}_2 \times \mathbb{Z}_2$ as the symmetry group of the resulting framework, since the choices above are independent.

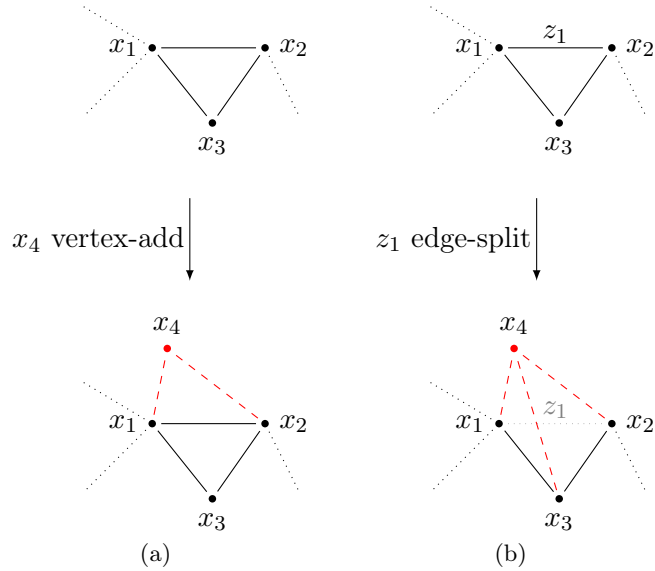

Figure 12: The two basic operations of the Henneberg sequence are illustrated.

In summary, constructing the 2-cycles can be done by two vertex-add operations and the choices with which one can perform these vertex-adds can be related to the addition of the groups $\mathbb{Z}_2$ as symmetry groups of the resulting frameworks. A natural question is thus whether one can find a general scheme to relate a symmetry group to the Henneberg sequences defining the framework. Such a characterization, by allowing to resort to group theoretic techniques and the known classification of discrete groups, would greatly deepen our understanding of both directed and undirected formation control.

## 5 Decentralized control with directed graphs

In this section, we formally introduce the decentralized control model used. We consider kinematic models of the form
$$\dot{x}_i = u_i$$

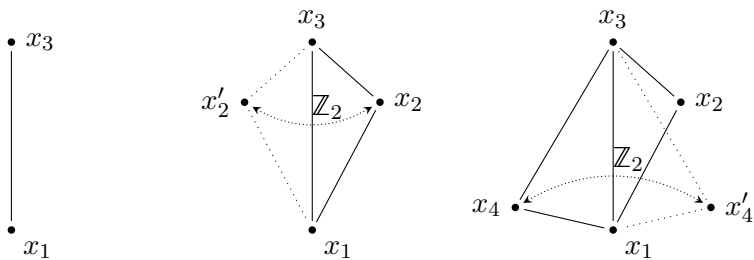

Figure 13: A Henneberg sequence for the 2-cycles. Starting with two nodes $(x_1, x_3)$ joined by an edge, we used two vertex-add operations to introduce $x_2$ and $x_4$. Each vertex-add operation adds a copy of the group $\mathbb{Z}_2$ as the symmetry groups of the formation. The two operations being commutative, the final symmetry group is $\mathbb{Z}_2 \times \mathbb{Z}_2$.

The formation control problem is then the following:
*Formation control problem*: Given a graph $G = (V, E)$, and a set of target distances $d \in \mathcal{L}$, find controls $u_i$, which respect the information flow described by $G$, such that for almost all initial conditions in $\mathbb{R}^{2n}$, the system stabilizes to a framework $p$ with $\delta(p)|_E = d$ — in other words, as $t \to \infty$, the inter-agent distances are given by $d$.

In order for the set of frameworks with inter-agent distances given by $d$ to be finite, the underlying graph has to be rigid. If it is moreover minimally rigid, no edges in the graph can be spared without losing rigidity.

## 5.1 Feasibility

Determining which frameworks are achievable by agents in a formation control problem can be a somewhat delicate problem when the underlying information flow is given by a *directed* graph. In particular, the fact that the graph is rigid, or even minimally rigid, is not enough to ensure that the agents will be able to reach the desired formation, and rigidity is not always needed. We will show some examples below that illustrate some of these issues and refer the reader to [BS03, HADB07] and references therein for further information.

Let us introduce some of the terminology used: an agent is called a *leader* if it is not following any other agent. An agent is a *coleader* if it follows some agents and is himself being followed by agents. An agent is a *follower* if it just follows agents and is not being followed.

Observe that any formation which has a follower with a single leader will not be rigid, but nevertheless may be feasible as is the case with the formation depicted in Figure 14a. The formation of Figure 14b, which has a follower with three coleaders is not feasible. Indeed, agent $x_3$ has no means to influence the position of agents $x_1, x_2$ and $x_4$ and these three will settle generically at a position which is not compatible with the constraints that $x_3$ has to satisfy. The problem with this formation is easily identifiable: $x_3$ is following more than 2 agents and the underlying graph is *minimally rigid*. In contrast, in the formation of Figure 14c, $x_3$ similarly has an outvalence greater than 2 but nevertheless can satisfy its three constraints since the graph without $x_3$ is rigid and the constraints are consistent—a consequence of the fact that the underlying graph is *redundantly rigid*.

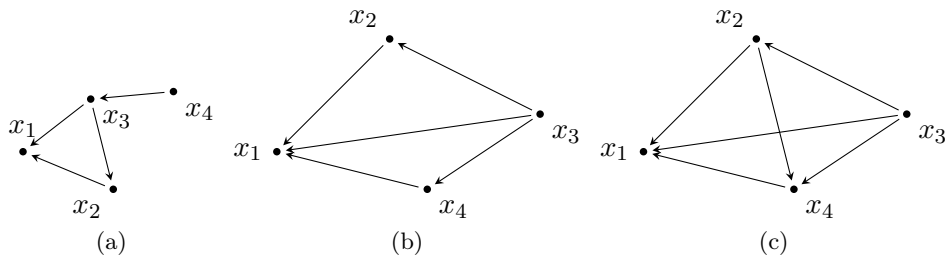

Figure 14

## 5.2 Control problem

### 5.2.1 Distributed formation control

We describe a general dynamical model for formation control with directed information flow. The model has to reflect the following three conditions:

1. Each agent is only aware of the target lengths it has to achieve.

2. Each agent is only aware of the position of its (co)leaders.

3. The agents do not share a common reference frame and only know the *relative* positions of their coleaders with respect to their own positions.

The first condition states that agents are only aware of the objective they have to achieve, and not the objective of the formation as a whole. In particular, agents do not know how many other agents are in the formation. The second condition says that the information flow in the system is given by the underlying graph. These first two conditions are at the core of decentralized control problems. The third condition enforce that the resulting system is invariant under the $SE(2)$ action described in Section 4. For example, if $x_i$ is such that $(i,j),(i,k) \in E$, then $u_i$ is allowed to depend on $x_j - x_i$ and $x_k - x_i$.

To the best of our knowledge, these assumptions are verified by most models that have appeared in the literature. An exception is the discussion in [YADF09], which separates the design stage from the dynamics stage and hence each agent is designed with a complete knowledge of the desired final formation. The first condition is then not satisfied. This has yielded only results of local nature.

### 5.2.2 The model

If edge $l$ connects agents $i$ and $j$, we write

$$e_l = \|x_i - x_j\|^2 - d_l,$$

the error in square edge length. We follow this convention, without loss of generality, in order to have the $d_i$ enter the dynamics linearly. We have the following model:

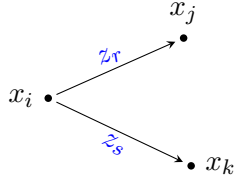

Figure 15

1. **Agent with outvalence of 1**: If agent $i$ has a unique leader $j$ and the target length for $\|x_i - x_j\|$ is $d_s$, then
$$\dot{x}_i = u(d_s; e_s)(x_j - x_i) \tag{8}$$

2. **Agent with outvalence of 2**: Assume that agent $i$ has a two leaders $j, k$ and that the target edge lengths are $d_r$ and $d_s$ respectively (see Figure 15). We will consider two cases, according to whether agent $i$ is able to measure explicitly how agent $k$ and agent $l$ are positioned relative to each other:

    (a) **Distance only:** In this case, the model is an extension to two variables of the model for an agent with a single leader:
    $$\dot{x}_i = u_1(d_r, d_s; e_r, e_s)(x_j - x_i) + u_2(d_r, d_s; e_r, e_s)(x_k - x_i) \tag{9}$$

    (b) **Distance and angle:** Define $z_s = x_j - x_i$ and $z_r = x_k - x_i$. The model is
    $$\dot{x}_i = u_1(d_r, d_s; e_r, e_s, z_s \cdot z_r)(x_j - x_i) + u_2(d_r, d_s; e_r, e_s, z_s \cdot z_r)(x_k - x_i). \tag{10}$$

    Knowing $d_r, d_s, e_r, e_s$ as well as the inner product $z_r \cdot z_s$ indeed allows agent $i$ to reconstruct the relative positions of agents $j$ and $k$ using simple trigonometric rules.

There are similar laws for agents with a higher outvalence. For the reasons mentioned, we only consider in the paper agents with an outvalence of at most 2. We always have the additional implicit condition on the $u_i$'s that the resulting differential equation admits a solution on a long enough time interval.

Equations (9) and (10) allow agents with two leaders to treat them differently. We will also consider here models that prevent this differentiation:

$$\dot{x}_i = u(e_r, e_s, z_r \cdot z_s)z_r + u(e_s, e_r, z_r \cdot z_s)z_s. \tag{11}$$

Indeed, observe that substituting $e_r$ for $e_s$ and $z_r$ for $z_s$ leaves equation (11) unchanged. We conclude this section by verifying that the model is invariant under the $SE(2)$ action introduced above.

**Lemma 1.** *The dynamical system described by Equations (8), (9) and (10) defines a smooth dynamical system on $E^n$.*

*Proof.* Since $E^n$ is a quotient space, it is sufficient to check that the dynamics is equivariant under the group action, or in other words that if $\dot{x} = F(x)$ and $S \in SE(2)$, then

$$S\dot{x} = F(Sx) = SF(x).$$

The result is a simple consequence of the fact that the $e_i$ are clearly invariant under an action of $S$, and so are the inner products $z_i \cdot z_j$. ∎

The description of the system in terms of the $z$ variables is redundant. Because the LS-category of the state space of the system is strictly greater than 1, we cannot write a unique set of differential equations in terms of non-redundant variables, such as the edge lengths, to describe the evolution of the system.

We show that the assumptions on the information flow put constraints on the control terms. Given a rigid graph and $d \in \mathcal{L}$, we call *design equilibria* the frameworks such that $e_i = 0$.

**Theorem 5.** *The design equilibria of a leaderless infinitesimally rigid formation are such that $u_i(0, 0, w) = u_j(0) = 0$.*

The theorem states that the system cannot be at a design equilibrium while undergoing a rigid transformation; e.g. the system cannot be such that $e_i = 0$ while the formation undergoes a translation. We need two preliminary lemmas. We first derive a useful formula for the dynamics in terms of the $z$ variables. Let

$$z = [z_1^T, z_2^T, \ldots, z_M^T]^T.$$

**Lemma 2.** *Given a graph $G$ with edge-adjacency matrix $A_e$, define the diagonal matrix $D$ with*
$$D_{ii} = u(e_i)$$
*if edge $i$ originates from an agent with a single leader and*
$$D_{ii} = u_i(e_i, e_j, z_i^T z_j)$$
*if edge $i$ originates from an agent with two leaders and $z_j$ links the origin of $i$ to the other leader. We have*
$$\dot{z} = A_e^{(2)} D^{(2)} z. \tag{12}$$

We recall that $A_e^{(2)}$ is a shorthand notation for $A_e \otimes I$ where $I_2$ is the $2 \times 2$ identity matrix and $\otimes$ the Kronecker tensor product.

*Proof.* Notice that a cycle in the underlying graph of a framework yields a linear relation that the $z$ variables have to satisfy. Hence, the $z$ variables are not independent if the underlying graph is not a tree and there are many non-equivalent ways of writing the dynamics in terms of these variables. We will verify formula (12) row by row. Consider the following generic situation depicted in the figure below: $x_1$ follows $x_2$ and $x_3$, $x_2$ follows $x_4$ and $x_5$ and $x_3$ follows $x_6$, with $z_1 = x_2 - x_1$ and $x_2 = x_3 - x_2$. We have

$$\begin{aligned}
\dot{z}_1 &= \dot{x}_2 - \dot{x}_1 \\
&= u_1(e_3, e_4, z_3^T z_4) z_3 + u_2(e_3, e_4, z_3^T z_4) z_4 - u_3(e_1, e_2, z_1^T z_2) z_1 - u_4(e_1, e_2, z_1^T z_2) z_2.
\end{aligned}$$

Observe that the first row of the edge adjacency matrix corresponding to this formation is $\begin{bmatrix} -1 & -1 & 1 & 1 & 0 & \ldots & 0 \end{bmatrix}$. Hence, this equation is indeed the first row of (12).

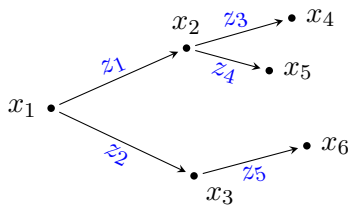

We have similar relations for the other rows of (12). ∎

**Lemma 3.** *Let $G = (V, E)$ represent a leaderless formation with $|V| = n$ and $|E| = m$. Then there exists a matrix $K \in \mathbb{R}^n \times \mathbb{R}^m$ of rank $n$ such that*
$$A_e = A_m K.$$

*Proof.* Fix a vertex $s$ in $V$ and let $j$ be an edge originating from this vertex. Since $G$ is leaderless, we know that such a $j$ exists. From the definition of $A_e$, we have that $A_{e,kj}$ is 1 if edge $k$ starts at $s$, $-1$ if it ends at $s$ and zero otherwise. Observe that, from the definition of $A_m$, we have the same relations for $A_{m,ks}$: $A_{m,ks}$ is 1 if edge $k$ starts at vertex $s$, $-1$ if it ends at vertex $s$ and zero otherwise. Hence if edge $j$ starts at vertex $s$, the corresponding columns in $A_e$ and $A_m$ are similar. Because the above is true for any vertex $s \in V$, and because every vertex has at least one edge originating from it, we have proved that $K$ exists and it contains the $n \times n$ identity matrix in its column span, which yields the result. ∎

We write $u_x$ for $\frac{\partial}{\partial x} u(x, y, z)$ where $x, y, z$ are dummy real variables. We can now prove Theorem 5:

*Proof of Theorem 5.* Let us write $\mathrm{diag}(D)$ for the vector whose entries are the diagonal entries of $D$, where $D$ is as defined in Lemma 2. Define the $m \times 2n$ matrix $Z$ by

$$Z = \begin{bmatrix} z_1 & 0 & \ldots & 0 \\ 0 & z_2 & \ldots & 0 \\ 0 & \ldots & \ddots & \vdots \\ 0 & 0 & \ldots & z_m \end{bmatrix}^T. \tag{13}$$

At a design equilibrium, we have $e_i \equiv 0$ by definition and thus the entries of $D$ are $u_i(0, 0, w)$ or $u_j(0)$ depending on the number of coleaders of the agent. Since $e_i \equiv 0$, we have $\frac{d}{dt}(z_i \cdot z_i) = 0$. A short computation using Equation (12) yields

$$\frac{d}{dt} \begin{bmatrix} z_1 \cdot z_1 \\ z_2 \cdot z_2 \\ \vdots \\ z_m \cdot z_m \end{bmatrix} = ZA_e Dz$$

$$= ZA_e Z^T \mathrm{diag}(D) = 0$$

Because the formation is leaderless, we have, using Lemma 3, $A_e = A_m K$ with $K$ of full rank. The rigidity matrix of the framework can be written as

$$R = ZA_m^{(2)}$$

where we recall that $A_m$ is the mixed-adjacency matrix of the underlying graph. We thus have

$$ZA_m KZ^T \mathrm{diag}(D) = RKZ^T \mathrm{diag}(D) = 0.$$

When the formation is infinitesimally rigid and leaderless, $R$, $K$ and $Z$ are of full rank. Hence, $RK \in \mathbb{R}^{m \times 2m}$ is of full row rank while $Z^T \in \mathbb{R}^{2m \times m}$ is of full column rank. Hence, $RKZ^T$ is of full rank and we conclude that $\mathrm{diag}(D) = 0$. ∎

The first condition states that agents are only aware of the target distances to their coleaders, and that they are unaware of the rest of the formation. Hence, agents do not know what the relative positions of their coleaders with respect to each other is at an equilibrium. We elaborate now on the implications of this fact on the model.

**Lemma 4.** *The 2-cycles formation is infinitesimally rigid for edge lengths $d \in \mathcal{L}_0$.*

*Proof.* The rigidity matrix of the 2-cycles is given by:

$$R = \begin{bmatrix} z_1 & -z_1 & 0 & 0 \\ 0 & z_2 & -z_2 & 0 \\ -z_3 & 0 & z_3 & 0 \\ 0 & 0 & z_4 & -z_4 \\ z_5 & 0 & 0 & -z_5 \end{bmatrix}^T \tag{14}$$

where the $z_i$ are defined by Equation (6). Recall that the rigidity matrix $R$ of a framework can be expressed as
$$R = Z A_m^{(2)}.$$
where we recall that $A_m$ is the mixed adjacency matrix, $A_m^{(2)} = A_m \otimes I$, $I$ being the $2 \times 2$ identity matrix and $Z$ is given in equation 13.

In the case of the 2-cycles, the mixed adjacency matrix $A_m \in \mathbb{R}^{5 \times 4}$ is of rank 3. The cokernel[3] of $A_m$ is spanned by $[0,0,1,1,1]^T$ and $[1,1,1,0,0]$. Hence, the cokernel of $A_m^{(2)}$ is four dimensional and spanned by the vectors $[0,0,1,1,1]^T \otimes [1,0]^T$, $[0,0,1,1,1]^T \otimes [0,1]^T$, $[1,1,1,0,0]^T \otimes [1,0]^T$ and $[1,1,1,0,0]^T \otimes [0,1]^T$

The matrix $Z$ is of full rank unless $z_i = 0$ for some $i$, which corresponds to two agents superposed. We thus have that $Z$ is of full rank for formations in $\mathcal{L}_0$. The kernel of $Z$ is given by the relations $z_1 + z_2 + z_3 = 0$ and $z_3 + z_4 + z_5 = 0$. Because $Z$ is of full row rank, $R$ is of full (row) rank if $A_m^{(2)}$ maps *onto* the coimage of $Z$. It is readily verified to be the case from the above relations describing the cokernel of $A_m^{(2)}$ and the kernel of $Z$. ∎

**Proposition 3.** *Let $d \in \mathcal{L}_0$. The feedback functions $u_i(d_i, d_j; x, y, z)$ have to be such that $\frac{\partial u_i}{\partial z}(d_i; d_j; 0, 0, z) = 0$ for all $z$ in the domain of $u_i$.*

We can informally understand this proposition as a consequence of the fact that an agent with two leaders does not know what the angle between its two leaders is at a design equilibrium, since this angle is a function of information the agent does not have access to.

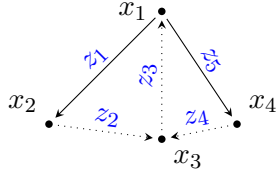

Figure 16

*Proof.* Let us assume that agent 1 has two coleaders. Without loss of generality, we consider the situation of Figure 16, as agent 1's dynamics does not depend directly on agents beyond its coleaders. The formation depicted is the 2-cycles.

The formation is at a design equilibrium with edge lengths $\|z_i\|^2 = d_i$. The differential equations for $z_1$ and $z_5$ are

$$\begin{cases} \dot{z}_1 &= u_2(e_2)z_2 - u_1(e_1, e_5, z_1 \cdot z_5)z_1 - u_5(e_1, e_5, z_1 \cdot z_5)z_5 \\ \dot{z}_5 &= u_4(e_4)z_4 - u_1(e_1, e_5, z_1 \cdot z_5)z_1 - u_5(e_1, e_5, z_1 \cdot z_5)z_5 \end{cases}$$

Let us consider a one-parameter family in $\mathcal{L}_0$ with the property that the associated frameworks all have $d_1$ and $d_5$ fixed to a constant value while the angle between $z_1$ and $z_5$ varies. Such a family exists because $d$ is in the interior of $\mathcal{L}$. Varying $d_3$ while keeping the other $d_i$'s constant yields such a family with the notation of Figure 16. We parametrize this family with $s \in (-\varepsilon, \varepsilon)$, $\varepsilon > 0$ and denote by $\gamma(s)$ the angle between $z_1$ and $z_5$.

Because the 2-cycles is infinitesimally rigid when $d \in \mathcal{L}_0$ by Lemma 4, we can use Theorem 5 to deduce that

$$u_1(d_1, d_5; 0, 0, \gamma(s)d_1 d_5) = u_5(d_1, d_5; 0, 0, \gamma(s)d_1 d_5) = 0$$

for almost all $s \in (-\varepsilon, \varepsilon)$. Hence, $u_{1,z} = u_{5,z} = 0$. The same argument can be repeated for every agent with two coleaders. ∎

We conclude by stating the *compatibility* conditions that $u_i$ have to satisfy in order to define a valid formation control system:

**Definition 7.** *An set of feedback control laws $u_i$ is* compatible *with the formation control problem if*

1. *$u_i(d_j; e_j)$ is such that $u_i(d_j; 0) = 0$ if agent $i$ has one co-leader.*

2. *$u_i(d_j, d_k; e_j, e_k, z_j^T z_k)$ is such that $u_i(d_j, d_k; 0, 0, z) = 0$ for all $z$ if agent $i$ has two co-leaders.*

---
[3]The cokernel of a linear map $f : A \to B$ is the quotient space $B/\operatorname{im}(f)$. Its coimage is $A/\ker(f)$.

## 5.3 Linearization of the dynamics

Given a formation control problem, the equations in $z$ variables are clearly redundant, as the quantity $\sum z_i$ is zero along any cycle in the underlying graph. They have the advantage, however, of being invariant under translation, which renders some proofs below more transparent. In this section, we will look at the linearization of the dynamics in the these variables. We let $F(z)$ be the right-hand side of the differential equation describing the system:
$$\dot z = F(z).$$

**Proposition 4.** *Set*
$$z'_i = (u_{1x} z_i + u_{2x} z_j)$$
*and*
$$z'_j = (u_{1y} z_i + u_{2y} z_j)$$
*if $z_i$ originates from an agent with two co-leaders given by $z_i$ and $z_j$, and*
$$z'_i = u_x z_i$$
*if $z_i$ originates from an agent with a single co-leader. Define*
$$Z' = \begin{bmatrix} z'_1 & 0 & 0 & \cdots & 0 \\ 0 & z'_2 & 0 & \cdots & 0 \\ 0 & \cdots & & \ddots & \vdots \\ 0 & 0 & \cdots & 0 & z'_m \end{bmatrix}. \tag{15}$$

*The Jacobian at a design equilibrium of a formation control system is*
$$\frac{\partial F}{\partial z} = A_e^{(2)} Z'^T Z. \tag{16}$$

*Proof.* We first observe that
$$\frac{\partial}{\partial z_i} u_1(e_i, e_j, z_i^T z_j) z_i = 2 u_{1x} z_i z_i^T + u I_2 + 2 u_{1z} z_i z_j^T$$
$$\frac{\partial}{\partial z_j} u_1(e_i, e_j, z_i^T z_j) z_i = 2 u_{1y} z_i z_j^T + 2 u_{1z} z_i z_i^T$$
$$\frac{\partial}{\partial z_i} u_2(e_i, e_j, z_i^T z_j) z_j = 2 u_{2x} z_j z_i^T + 2 u_{2z} z_j z_j^T$$
$$\frac{\partial}{\partial z_j} u_2(e_i, e_j, z_i^T z_j) z_j = 2 u_{2y} z_j z_j^T + 2 u_{2z} z_j z_j^T$$

where we omitted the arguments of the functions on the right-hand side in order to keep the notation simple. Recall from Proposition 3 that at a design equilibrium $u_z = 0$. Hence, if $z_i$ originates from an agent with two co-leaders with

$$\dot{z}_i = F_i = \ldots + u_1(e_i, e_j, z_i^T z_j) z_i + u_j(e_i, e_j, z_i^T z_j) z_j$$

then:

$$\begin{aligned}\frac{\partial F_i}{\partial z_i} &= -2u_{1x} z_i z_i^T - 2u_{2x} z_j z_i^T \\ &= z_i' z_i^T.\end{aligned}$$

If $z_j$ originates from an agent with a single leader, we have:

$$\frac{\partial F_j}{\partial z_j} = u_x z_j z_j^T = z_j' z_j^T.$$

Putting the equations above together, we get the result after some simple algebraic manipulations. ∎

Recall that the problem is invariant under an action of the Euclidean group $SE(2)$ on $\mathbb{R}^2$ and that, in addition, the $z$ variables are not independent. As a consequence, the Jacobian of the system expressed in the $z$ variables at an equilibrium will always have multiple zero eigenvalues. This degeneracy of the Jacobian does not reveal anything about the dynamics of the system beyond the fact that equilibria are part of manifold of equilibria. For all practical purposes, this connected set of equilibria can be thought of as a single equilibrium by taking the quotient by the action of the group. Under this quotient, this degeneracy of the Jacobian disappears.

To address this issue, the following result gives a computational tool to evaluate the eigenvalues of the Jacobian that do not correspond to the action of the Euclidean group while working in the more convenient inter-agent distance coordinates. Observe that in the Proposition below $J \in \mathbb{R}^{m \times m}$.

**Corollary 1.** *Let $G$ be the graph of a minimally rigid formation. The eigenvalues of the Jacobian at a design equilibrium are the eigenvalue zero with algebraic multiplicity $2n+3-m$ and the eigenvalues of*

$$J = Z A_e Z'^T$$

*Proof.* The result is a consequence of Theorem 1.3.20 in [HJ90] when applied to Proposition 4. ∎

# 6  Conclusion

We have investigated some geometric properties of formation control. The main contributions of this paper can be summarized as follows:

- We have described the geometry of the space $E^n$ of configurations of $n$ agents in the plane and shown that this space was $E^n = \mathbb{C}\mathrm{P}(n-2) \times \mathbb{R}^+$

- We have shown how a global topological characteristic of $E^n$, namely the LS-category, relates to the coordinates used by the agents in the ensemble, called localized coordinates. By doing so, we have derived a lower bound on the number of non-congruent frameworks with similar edge lengths when $\mathcal{L}$ is null-homotopic.

- We have established a number of conditions that a feedback control law needs to satisfy in order to model a formation control problem appropriately.

In particular, formations with 3 and 4 agents have been investigated in details as examples of the main results. We also conjectured the existence of a general procedure relating the Henneberg sequence to discrete symmetry groups acting on frameworks.

The results of this paper will be used extensively in part II, whose focus is on dynamical properties of formation control. We develop there a framework to analyze the range of behaviors achievable by the system while respecting the information flow of the underlying graph. We then apply this framework to show that there are no compatible feedbacks, as described in Definition 7, that satisfactorily stabilize the 2-cycles.